\documentclass [12] {article}

\usepackage[english]{babel}
\usepackage{amsthm}
\usepackage{epsfig, graphicx, xcolor}
\usepackage{amssymb}
\usepackage{amsmath, amsfonts}

\addtocounter{page}{0}

\makeatletter
\renewcommand\@biblabel[1]{#1.}
\makeatother

\begin{document}

\title{
Optimal recovery of convolutions of $n$ functions\\ according to linear information}

%

\author{V. F.~Babenko, M. S.~Gunko}

\date{Oles Honchar Dnepropetrovsk National University}

\maketitle

\begin{abstract}

We found the optimal linear information and the optimal method of
its use to recovery of the convolution of n functions on some classes of $2\pi$ - periodic functions.

\medskip


\end{abstract}

{\large Let $C$ and $L_p$, ${1\le p\le{\infty}},$ be the spaces $2\pi$ - periodic functions 
${f: \mathbb{R}\to \mathbb{R}}$ endowed with the corresponding norms $\|\cdot\|_C$  and $\|\cdot\|_{Lp}$.
Let $M_1,...,M_n \subset L_1$ be some classes of functions; $x_1 \in M_1$, $x_2 \in M_2$, ..., $x_n \in M_n$;
$$
(x_1* x_2)(\tau)=\int\limits_0^{2\pi} x_1(\tau - t)x_2(t) dt
$$
be the convolution of two functions, $x_1$ and $ x_2$, and
$$
(x_1* x_2*...*x_n)(\tau)=(x_1* (x_2*...*x_n))(\tau)
=
$$
$$
=\int\limits_0^{2 \pi} \cdot \cdot \cdot \int\limits_0^{2 \pi}
x_1(\tau - t_1-...-t_{n-1}) x_2(t_1)...x_n(t_{n-1}) dt_1 dt_2 ... dt_{n-1}
$$
be the convolution of $n$ functions $x_1,...,x_n$.

Suppose that for $j=1,...,n$ there are given sets
$ T_j=(T_{j,1},...,T_{j,m_j})$ of linear continuous functionals
$$T_{j,l}: {\rm span}(M_j) \to \ \mathbb{R},\;\;  l=1,...,m_j.$$

Vectors
$$
T_j(x_j)=(T_{j,1}(x_j),...,T_{j,m_j}(x_j))\in \mathbb{R}^{m_j},\;x_j \in M_j,\; j=1,...,n,
$$
will be called the linear information about $x_1$, $x_2$,..., $x_n$ of the type $(m_1,...,m_n)$ (or $(m_1,...,m_n)$-information).
An arbitrary mapping
\[\Phi:{\mathbb{R}^{m_1}
\times\ldots\times \mathbb{R}^{m_n}}\rightarrow
L_p\]
 we will be called a method of the recovery of convolution ${x_1* x_2*...*x_n}$  with respect to the $(m_1, ..., m_n)$-information in the $L_p$ space.

 We set:
$$
R(x_1,...x_n; T_1, ... , T_n;  \Phi)(t)= (x_1* x_2*...*x_n)(t) - \Phi(T_1(x_1),..., T_n(x_n))(t),
$$
\begin{equation}
 \label{eqn1}
R(M_1, ..., M_n; T_1, ..., T_n; \Phi;L_p)
=
\sup_{x_j \in M_j,\atop
j = 1,..., n} \| R(x_1,...,x_n; T_1, ..., T_n;\Phi)(\cdot) \|_{L_p},
\end{equation}
\begin{equation}
 \label{eqn2}
R(M_1, ..., M_n; T_1, ..., T_n;L_p)=\inf_\Phi
R(M_1, ..., M_n; T_1, ..., T_n; \Phi; L_p),
\end{equation}
\begin{equation}
 \label{eqn3}
R_{m_1, ..., m_n} (M_1,...,M_n;L_p)=\inf_{T_1, ..., T_n}
R(M_1,...,M_n;T_1,..., T_n;L_p)
\end{equation}
($\inf\limits_{\Phi}$ is taken over all possible methods of recovery, and $\inf\limits_{T_1, ..., T_n}$ is taken over all possible sets of functionals that give
$(m_1, ..., m_n)$-- information about $(x_1, ..., x_n)$;
\begin{equation}
 \label{eqn4}
R_{N}(M_1,...,M_n;L_p)=\inf_{m_1+...+m_n=N}
R_{m_1,...,m_n}(M_1,...,M_n;L_p).
\end{equation}
Quantity \eqref{eqn1} will be called the error of the method $\Phi$ of recovery of the convolution of $n$ functions ${x_1* x_2*...*x_n}$
on classes $M_1 , ...,M_n$ according to information $T_1 , ...,T_n$ in the space $L_p$;
quantity \eqref{eqn2} will be the called the optimal error of recovery of ${x_1* x_2*...*x_n}$ on classes $M_1 , ...,M_n$
according to given information of the type $(m_1,...,m_n)$;
quantity \eqref{eqn3} will be called the optimal error of recovery of ${x_1* x_2*...*x_n}$ on  $M_1, ..., M_n$
with respect to the $(m_1,...,m_n)$ - information; and quantity  \eqref{eqn4} will be called the optimal error of recovery of
${x_1* x_2*...*x_n}$ on  $M_1,...,M_n$ according to the information of the volume $N$.

If, for given $T_1, ..., T_n$, a method $\Phi^*$ that realizes $\inf\limits_{\Phi }$ on the right-hand side of \eqref{eqn2} exists,
then $\Phi^*$  will be called the optimal method of using this information.
If sets of functionals $T^*_1,...,T^*_n$ that realize lower bound in the right-hand side of \eqref{eqn3} exist, then we will call them the optimal $(m_1, ...,m_n)$ - information for recovery of ${x_1* x_2*...*x_n}$ on  $M_1,...,M_n$ in the space $L_p$.

Numbers $m_1^0,...,m_n^0$, for which the infimum is realized in \eqref{eqn4}, will be called the volumes of information about $x_1,...,x_n$,
at the same time, the optimal $(m_1^0,...,m_n^0)$ - information will be called the optimal information about $x_1,...,x_n$ of the volume $N$.

We will study the next problem of optimal recovery of  convolution of $n\;$ functions according to linear information about the functions.

Let the classes $M_1,\ldots ,M_n$ and the numbers $p\in [1,\infty ]$, $N\in \mathbb{N}$ be given.
Find the value \eqref{eqn4}, optimal volumes $m^*_1,\ldots ,m^*_n$ of information ($m^*_1+\ldots +m^*_n=N$), optimal $(m^*_1,\ldots ,m^*_n)$-information  $T^*_1,...,T^*_n$, and optimal method $\Phi^*$ of its use.

The problem of optimal recovery of convolution of two functions from various classes was studied in  {\cite{vf1989}}. The first results on investigation of this problem were also obtained in this paper.
For the other results in this direction, see {\cite{vfr1991}}.

We define the sets
$M_j(T_j)$ as follows:
$$
M_j(T_j)=\left\{ x_j \in M_j: T_j(x_j) = 0 \right\},j=1,...,n.
$$
Let
$$
M(T_j):=M_1 \times... \times M_{j-1} \times M_j(T_j)
\times M_{j+1} \times... \times M_n, j=1,...,n .
$$

The following lemma gives a lower estimate for the value  \eqref{eqn1}, and, hence, for the value  \eqref{eqn2}.

{\bf Lemma 1.}
\label{Lemma 1}
{\it Let classes $M_1,\ldots,M_n$ be convex and centrally - symmetric.
Then, for any sets of functionals $T_1,..., T_n$ and for any method of
 recovery $\Phi$
$$
R(M_1,...,M_n;T_1,..., T_n;\Phi;L_p ) \ge
$$
$$
\ge \max_{
j=1,...,n
}
\sup_{
(x_1,..., x_n) \in M(T_j)
}
\| (x_1* x_2*...*x_n)(\cdot)\|_{L_p}.
$$
}

$\Box$
Let us prove that
$$
R(M_1,...,M_n;T_1,.., T_n;\Phi;L_p)
\ge
\sup_{
(x_1,...,x_n) \in M(T_1)
}
\| (x_1*...*x_n)(\cdot)\|_{L_p}.
$$
We have (below, $\theta$ is null-element of the space $\mathbb{R}^{m_1}$)
$$
R(M_1,...,M_n;T_1,.., T_n;\Phi;L_p)
\ge
$$
$$
\ge\sup_{
(x_1,...,x_n) \in M(T_1)
} \|(x_1* x_2*...*x_n)(\cdot)-\Phi(\theta, T_2(x_2),...,T_{n}(x_{n}))(\cdot) \|_{L_p}\ge
$$
$$
\ge
\sup_{
(x_1,...,x_n) \in M(T_1)
}
\max \left\{
\|(x_1* x_2*...*x_n)(\cdot)-\Phi(\theta, T_2(x_2),...,T_{n}(x_{n}))(\cdot) \|_{L_p},\right.
$$
$$
\left. \| -(x_1* x_2*...*x_n)(\cdot)-\Phi(\theta, T_2(x_2),...,T_{n}(x_{n}))(\cdot) \|_{L_p}\right\} \ge
$$
$$
\ge\sup_{
(x_1,..., x_n) \in M(T_1)
}
\| (x_1* x_2*...*x_n)(\cdot)\|_{L_p}.
$$
Similarly, we obtain that for $j =2,..., n $
$$
R(M_1,...,M_n;T_1,.., T_n;\Phi;L_p)
\ge
\sup_{
(x_1,...,x_n) \in M(T_j)
}
\| (x_1* x_2*...*x_n)(\cdot)\|_{L_p}.
$$
Lemma is proved.
$\Box$

Let $F_p$ be the unit ball in $L_p$, $K \in L_1$, $\int\limits_0^{2 \pi} K dt \ne 0$.
Let us denote by $K*F_p$ the class of functions  of the form ${x= K*\psi}$,
$\psi \in F_p$.
 Further we will suppose that $M_j =K_j*{{F_p}_j}$, $ j=1,...,n $  , where $K_j\in L_1$.

Denote by  $d_N(M,C)$ the Kolmogorov $n$-width of the set $M$ in the space $C$ (see, e. g., {\cite {knp1984}},
 ё.109).

{\bf Theorem 1.}
\label{Th1}
{\it Let $K_1,...,K_n \in L_1$.
Then, for any sets of linear functionals $T_j=(T_{j,1},...,T_{j,m_j})$, $j=1,...,n,$ $T_{j,l}:{\rm span}\,(K_j*F_1)\to \mathbb{R}$, $l=1,...,m_j$
 and for any method of recovery $\Phi$
$$
 R(K_1*F_1, ..., K_n*F_1; T_1, ..., T_n; \Phi;L_1)
\ge
$$
$$\ge
d_{\min\{m_1,\ldots,m_n\}}((K_1\ast K_2\ast...\ast K_n)(-\cdot)\ast F_{\infty}, {C}).
$$
}
$\Box$ Let, for definiteness, $\min\{m_1,\ldots,m_n\}=m_1$. Due to the Lemma 1 we only need to prove the following inequality
$$
A_1:=\sup_{x_1\in (K_1*F_1)(T_1)\atop x_2\in K_2*F_1,..., x_n\in K_n*F_1 }\left\| (x_1\ast x_2\ast...\ast x_n)(\cdot)\right\|_{L_1}\ge d_{m_1}(K_1\ast K_2\ast...\ast K_n\ast F_{\infty}, {C}).
$$
We have
$$
A_1=\sup_{\psi_j\in F_1,\, j=1,...,n,\atop T_1(K*\psi_1)=\theta}\| (K_1\ast\psi_1)\ast(K_2\ast\psi_2)\ast...\ast(K_n\ast\psi_n)(\cdot)\|_{L_1}=
$$
$$
=\sup_{\psi_j\in F_1,\, j=1,...,n,\atop T_1(K*\psi_1)=\theta}\| (K_1\ast\ldots\ast K_n\ast\psi_1\ast...\ast\psi_n)(\cdot)\|_{L_1}=
$$
$$
=\sup_{\psi_j\in F_1,\, j=1,...,n,\atop T_1(K*\psi_1)=\theta}\sup_{\phi\in F_\infty}\int\limits_0^{2\pi} (K_1\ast\ldots\ast K_n\ast\psi_1\ast...\ast\psi_n)(t)\phi(t)dt=
$$
$$
=\sup_{\phi\in F_\infty}\sup_{\psi_j\in F_1,\, j=1,...,n,\atop T_1(K*\psi_1)=\theta}\int\limits_0^{2\pi} ((K_1\ast\ldots\ast K_n)(-\cdot)\ast \phi)(u)(\psi_1\ast...\ast\psi_n)(u)du=
$$
$$
=\sup_{\phi\in F_\infty}\sup_{\psi_j\in F_1,\, j=1,...,n,\atop T_1(K*\psi_1)=\theta}\max\limits_t\int\limits_0^{2\pi} ((K_1\ast\ldots\ast K_n)(-\cdot)\ast \phi)(u)(\psi_1\ast...\ast\psi_{n-1})(u-t)du=
$$
$$
=\sup_{\phi\in F_\infty}\sup_{\psi_j\in F_1,\, j=1,...,n-1,\atop T_1(K*\psi_1)=\theta}\int\limits_0^{2\pi} ((K_1\ast\ldots\ast K_n)(-\cdot)\ast \phi)(u)(\psi_1\ast...\ast\psi_{n-1})(u)du=\ldots=
$$
$$
=\sup_{\phi\in F_\infty}\sup_{\psi_j\in F_1,\, j=1,2,\atop T_1(K*\psi_1)=\theta}\int\limits_0^{2\pi} ((K_1\ast\ldots\ast K_n)(-\cdot)\ast \phi)(u)(\psi_1\ast\psi_2)(u)du=
$$
$$
=\sup_{\phi\in F_\infty}\sup_{\psi_j\in F_1,\atop T_1(K*\psi_1)=\theta}\int\limits_0^{2\pi} ((K_1\ast\ldots\ast K_n)(-\cdot)\ast \phi)(u)\psi_1(u)du.
$$

Functionals $T_{1,1},...,T_{1,{m_1}}$ functionals allow a continuous continuations $T'_{1,1},...,T'_{1,{m_1}}$ to the whole space $L_1$.
Let functionals $T'_{1,j}$ have the form
$$
T'_{1,j}(x_1)=\int\limits_0^{2\pi} x_1(t)g_{1,j}(t)dt, j=1,...,m_1,
$$
where $g_{1,j}$ are fixed functions from $L_\infty$.

The condition $T_1(K_1*\psi_1)=\theta$ means that for $j=1,\ldots,m_1$
$$
\int\limits_0^{2\pi}(K_1*\psi_1)(t)g_{1,j}(t)dt=0,
$$
or
$$
\int\limits_0^{2\pi}\int\limits_0^{2\pi}K_1(t-\tau)\psi_1(\tau)d\tau g_{1,j}(t)dt=0,
$$
or
$$
\int\limits_0^{2\pi}\psi_1(\tau)\int\limits_0^{2\pi}K_1(t-\tau)g_{1,j}(t)dtd\tau=0.
$$
Thus, the condition $T_1(K_1*\psi_1)=\theta$ means that $\psi_1\perp K_1(-\cdot)*g_{1,j}$ for any $j=1,\ldots ,m_1$.
 Therefore, taking into account the S. M. Nikol'sky's duality theorem (see {\cite{knp1984}, ё. 120} proposition 3.4.4), we obtain
$$
A_1=\sup_{\phi\in F_\infty}\sup_{\psi_1\in F_1,\atop \psi_1\perp K_1(-\cdot)*g_{1,j},\, j=1,\ldots ,m_1}\int\limits_0^{2\pi} ((K_1\ast\ldots\ast K_n)(-\cdot)\ast \phi)(u)\psi_1(u)du=
$$
$$
=\sup_{\phi\in F_\infty}E((K_1\ast\ldots\ast K_n)(-\cdot)\ast \phi); H(T_1))_C,
$$
where $E( x;H(T_1))_{C}$ is the best approximation of a function $x$ by subspace
$$
H(T_1)={\rm span} \{ K_1(-\cdot)*g_{1,1}, K_1(-\cdot)*g_{1,2},..., K_1(-\cdot)*g_{1,{m_1}}\}
$$
in the space $C$.

Since
$\dim H(T_1 )\le m_1$ we obtain taking into account the definition of Kolmogorov $n$-width, 
$$
A_1\ge  d_{m_1}(  (K_1*K_2*...\ast K_n)(-\cdot)\ast F_{\infty}, C).
$$
The Theorem 1 is proved.
$\Box$
Let us denote by $H^T_{2s-1}$, $s\in\mathbb{N}$ the set of trigonometric polynomials of the order not greater than $s-1$.
A kernel  $K$, which is continuous on $(0;2\pi)$ and is not a trigonometric polynomial, will be called a $CVD$ - kernel
(in short $K\in CVD$), if for every $\phi \in {C}$,
 $\nu(K*\phi)\le\nu(\phi)$, where $\nu(g)$ is the number of sing changes of the function $g$ on a period. Series of topics related to the theory of $CVD$ - kernels are presented in {\cite {msw1959}}, {\cite {ks1968}}.

Let $K\in CVD$ . Then it satisfies assumptions of theorem 4.1 from {\cite  {vf1987} } and, hence,
if $\varphi_st={\rm sign} \sin st$, $\sigma$ is the point of the absolute maximum or absolute minimum of the function
$K*\varphi_s$, then there exists  a unique polynomial \newline$P_{s,\sigma}=P_{s,\sigma}(K)\in {H^T}_{2s-1}$, interpolating $K(t)$ at points $\sigma+\frac{m\pi}{s},m\in \mathbb{Z}$, if all points $\sigma+\frac{m\pi}{s},m\in \mathbb{Z}$ are continuity at points of $K$. If $K$ has a discontinuity at null, and $0\in\{\sigma+\frac{m\pi}{s}|m\in \mathbb{Z}\}$,
then there exists a unique polynomial  $P_{s,\sigma}=P_{s,\sigma}(K)\in {H^T}_{2s-1}$,
interpolating $K$ at points $\sigma+\frac{m\pi}{s} \not \equiv 0 ({\rm mod} \; 2\pi)$.
The polynomial  $P_{s,\sigma}=P_{s,\sigma}(K)$ (see {\cite  {vf1987} } (theorems 2.4, 4.2 and \S
5), {\cite  {pa1979} }) is the best $L_1$ approximation  polynomial for $K$ and additionally
$$
\|K-P_{s,\sigma}\|_{L_1}=\|K(-\cdot)-P_{s,\sigma}(-\cdot)\|_{L_1}=\|K*\varphi_s\|_\infty=
$$
\begin{equation}
 \label{eqn8}
=d_{2s-1}(K*F_\infty;C)=d_{2s-1}(K(-\cdot)*F_\infty;C).
\end{equation}

Below, considering the problem of recovery of the convolution $x_1*\ldots *x_n$, where $x_1\in K_1*F_1,\ldots x_n\in K_n*F_1$, we will suppose that the kernels $K_1,\ldots, K_n$ are such that $K_1*\ldots *K_n\in CVD$.

Let $a_j(x)$, $b_j(x)$ be the Fourier coefficients of a function $x\in L_1$, i.e.
$$
a_j(x)=\frac{1}{\pi}\int\limits_0^{2\pi}x(t) \cos jtdt, \;\;\; b_j(x)=\frac{1}{\pi}\int\limits_0^{2\pi}x(t) \sin jt dt,
$$
$$
c_j(x)=\frac{a_j(x)-ib_j(x)}{2}, \;\;\; c_{-j}(x)=\frac{a_j(x)+ib_j(x)}{2},\;\;\; j=0,1,... .
$$
Set
$$\alpha_j=\frac{c_j(P_{n,\sigma}(K_1*K_2*...*K_n))}{(2\pi)^{n-1}c_j(K_1)c_j(K_2)...c_j(K_n)}, \;\;\;j=0,\pm1,...,\pm( s-1).
$$
Note that it follows from the assumption $K_1*\ldots *K_n\in CVD$  that all coefficients $c_j(K_1),c_j(K_2),...,c_j(K_n)$ are not equal to zero for every $j\in \mathbb{Z}$.

Let
\begin{equation}
 \label{eqn10}
{T_l}^*=(a_0(x_l),a_1(x_l),...,a_{N-1}(x_l),b_1(x_l),...,b_{N-1}(x_l)), l=1,...,n,
\end{equation}
and
\begin{equation}
 \label{eqn11}
\Phi^*({T_1}^*(x_1),...,{T_n}^*(x_n))(t)=\sum_{j=-(s-1)}^{s-1}\alpha_j c_j(x_1)...c_j(x_n)e^{ijt}.
\end{equation}
If $x_1=K*\psi_1\in  K_1*F_{1},...,x_n=K_n*\psi_n\in  K_n*F_{1}$, then  for $g_1,...,g_n \in L_1$
$$
c_j(g_1*g_2*...*g_n)=(2\pi)^{n-1}c_j(g_1)c_j(g_2)...c_j(g_n), j\in\mathbb{Z}.
$$
 Analogously to the transformation of the value $A_1$ in the proof of the 
Theorem 1, we obtain
$$
\|x_1*x_2*...*x_n(\cdot)-\Phi^*({T_1}^*(x_1),...,{T_n}^*(x_n);\cdot)\|_{L_1}=
$$
$$
=\|K_1*...*K_n*\psi_1*...*\psi_n - P_{s,\sigma}(K_1*...*K_n)*\psi_1*...*\psi_n\|_{L_1}=
$$
$$=\sup_{\phi\in F_\infty}\int\limits_0^{2\pi} [(K_1\ast\ldots\ast K_n\ast\psi_1\ast...\ast\psi_n)(t)-P_{s,\sigma}(K_1*...*K_n)*\psi_1*...*\psi_n]\phi(t)dt=
$$
$$
=\sup_{\phi\in F_\infty}\int\limits_0^{2\pi} [(K_1\ast\ldots\ast K_n)(-\cdot)(t)-(P_{s,\sigma}(K_1*...*K_n))(-\cdot)(t)]*\phi (t)\cdot(\psi_1*...*\psi_n)(t)dt.
$$
It follows that
$$
\|x_1*x_2*...*x_n(\cdot)-\Phi^*({T_1}^*(x_1),...,{T_n}^*(x_n);\cdot)\|_{L_1}\le
$$
$$
\le\sup_{\phi\in F_\infty}\| [(K_1\ast\ldots\ast K_n)(-\cdot)(t)-(P_{s,\sigma}(K_1*...*K_n))(-\cdot)]*\phi\|_C\|\psi_1*...*\psi_n\|_{L_1}\le
$$
$$
\le\sup_{\phi\in F_\infty}\| [(K_1\ast\ldots\ast K_n)(-\cdot)(t)-(P_{s,\sigma}(K_1*...*K_n))(-\cdot)]*\phi\|_C\le
$$
$$
\le\| (K_1*...*K_n)(-\cdot) - (P_{s,\sigma}(K_1*...*K_n))(-\cdot)\|_{L_1}\le\| K_1*...*K_n- P_{s,\sigma}(K_1*...*K_n)\|_{L_1}.
$$
Taking into account the relations~\eqref{eqn8}, we have
$$
\|x_1*x_2*...*x_n(\cdot)-\Phi^*({T_1}^*(x_1),...,{T_n}^*(x_n);\cdot)\|_{L_1}\le
$$
$$
\le d_{2s-1}(K_1*...*K_n*F_\infty;C).
$$

Thus,
$$
R(K_1\ast{{F_{1}}},K_2\ast{{F_1}},..., K_n\ast{{F_1}}; T_1, ..., T_n; \Phi^*;L_1)\le
$$
$$
\le d_{2s-1}(K_1*...*K_n*F_\infty;C)=\| K_1*...*K_n*\varphi_s\|_C.
$$

Taking into account the lower estimate, which is given by Theorem 1, and the monotonicity of the widths $d_{m}(M, {C})$ as $m$ growths,
we see that the following theorem holds true.

{\bf Theorem 2.}
\label{Th2}
{\it Let kernels $K_1,K_2,...,K_n$ be such that
$K_1*K_2*...*K_n\in CVD$, let $s\in \mathbb{N}$ and let $N=n(2s-1)$.
Then
$$
R_{n(2s-1)}(K_1\ast{{F_{1}}},K_2\ast{{F_1}},..., K_n\ast{{F_1}}; L_1)=R_{2s-1,\ldots,2s-1}(K_1\ast{{F_{1}}},K_2\ast{{F_1}},..., K_n\ast{{F_1}}; L_1)=
$$
$$
=d_{2s-1}(K_1\ast K_2\ast...\ast K_n\ast F_\infty, {C})
= \|K_1\ast K_2\ast...\ast K_n\ast \varphi_s \|_{C}.
$$
Additionally, the optimal information is given by equality \eqref{eqn10}, and the corresponding optimal method is given by equality \eqref{eqn11}.}



\addcontentsline{toc}{chapter}{}

}


%
%
\end{document}